\documentclass[reqno,12pt]{amsart}

\newtheorem{Theorem}{Theorem}[section]
\newtheorem{Proposition}[Theorem]{Proposition}
\newtheorem{Corollary}[Theorem]{Corollary}
\theoremstyle{remark}
\newtheorem{Remark}[Theorem]{Remark}
\numberwithin{equation}{section}

\newcommand\bm{\mathbf m}
\newcommand\bn{\mathbf n}
\newcommand\bk{\mathbf k}
\newcommand\bl{\mathbf l}

\allowdisplaybreaks

\author[Michael Schlosser]{Michael Schlosser$^*$}
\address{Fakult\"at f\"ur Mathematik, Universit\"at Wien,
Nordbergstra{\ss}e 15, A-1090 Vienna, Austria}
\email{michael.schlosser@univie.ac.at}
\urladdr{http://www.mat.univie.ac.at/{\textasciitilde}schlosse}
\thanks{$^*$ Partly supported by FWF Austrian Science Fund
grants \hbox{P17563-N13}, and S9607 (the second is part
of the Austrian National Research Network
"Analytic Combinatorics and Probabilistic Number Theory").}
\date{July 5, 2006}


\title[A new multivariable ${}_6\psi_6$ summation]
{A new multivariable
${}_{\boldsymbol 6}\boldsymbol\psi_{\boldsymbol 6}$ summation formula}

\begin{document}

\begin{abstract}
By multidimensional matrix inversion, combined with an
$A_r$ extension of Jackson's ${}_8\phi_7$ summation
formula by Milne, a new multivariable ${}_8\phi_7$
summation is derived. By a polynomial argument this
${}_8\phi_7$ summation is transformed to
another multivariable ${}_8\phi_7$ summation
which, by taking a suitable limit, is reduced to
a new multivariable extension of the nonterminating
${}_6\phi_5$ summation. The latter is then extended,
by analytic continuation, to a new multivariable extension
of Bailey's very-well-poised ${}_6\psi_6$ summation formula.
\end{abstract}

\maketitle

\section{Introduction}\label{sec0}

Bailey's~\cite[Eq.~(4.7)]{Ba} very-well-poised ${}_6\psi_6$
summation formula,
\begin{multline}\label{66gl}
{}_6\psi_6\!\left[\begin{matrix}q\sqrt{a},-q\sqrt{a},b,c,d,e\\
\sqrt{a},-\sqrt{a},aq/b,aq/c,aq/d,aq/e\end{matrix}\,;q,
\frac{a^2q}{bcde}\right]\\
=\frac {(q,aq,q/a,aq/bc,aq/bd,aq/be,aq/cd,aq/ce,aq/de)_{\infty}}
{(aq/b,aq/c,aq/d,aq/e,q/b,q/c,q/d,q/e,a^2q/bcde)_{\infty}},
\end{multline}
(see Section~\ref{secpre} for the notation)
where $|a^2q/bcde|<1$ (cf.\ \cite[Eq.~(5.3.1)]{GR}),
is one of the most important identities in special functions, with
applications to orthogonal polynomials, number theory, and combinatorics.

Several {\em multivariable} extensions of Bailey's formula exist
(all of them associated with various root systems),
including a couple of summations by Gustafson~\cite{G1,G2,G3},
a summation by van~Diejen~\cite{vD}, by the author~\cite{S3},
and by Ito~\cite{I1}.

In this paper we supply a new multivariable extension
of \eqref{66gl} to this list. We also provide
several other new summations. These include new multivariable
very-well-poised ${}_8\phi_7$ and ${}_6\phi_5$ summations.
The series obtained in this paper are of a slightly different
type than those usually labelled as $A_r$ series, but,
being closely related to these, we still decided to refer to them
as $A_r$ series, see Remark~\ref{rem}.

Our paper is organized as follows.
In Section~\ref{secpre} we introduce the notation and
review some common facts about basic hypergeometric series.
We first do this for the classical univariate case and then
for the multivariate case. In Section~\ref{secmmi}
we explain the concept of multidimensional matrix inversion
and list an explicit result, see Corollary~\ref{mmic}
which is needed in Section~\ref{secn87} to derive,
via multidimensional inverse relations and an $A_r$
${}_8\phi_7$ summation theorem by Milne~\cite{M1},
a new multivariable terminating very-well-poised
${}_8\phi_7$ summation.
A polynomial argument gives yet another multivariable
${}_8\phi_7$ summation. In Section~\ref{secn65}
we suitably specialize this summation in order to obtain
new multivariable very-well-poised ${}_6\phi_5$ summations.
Finally, in Section~\ref{secn66} we apply analytic
continuation (in particular, an iterated application of
Ismail's~\cite{AI,Is} argument) to deduce a new multivariable
extension of Bailey's very-well-poised ${}_6\psi_6$ summation.

\section{Preliminaries}\label{secpre}
\subsection{Notation and basic hypergeometric series}

Here we recall some standard notation for $q$-series,
and basic hypergeometric series (cf.\ \cite{GR}).

Let $q$ be a complex number such that $0<|q|<1$. We define the
{\em $q$-shifted factorial} for all integers $k$ by
\begin{equation}\label{ipr}
(a)_k:=\frac{(a)_{\infty}}{(aq^k)_{\infty}},
\end{equation}
where
\begin{equation*}
(a)_{\infty}:=\prod_{j=0}^{\infty}(1-aq^j).
\end{equation*}
For brevity, we occasionally employ for products the notation
\begin{equation*}
(a_1,\ldots,a_m)_k:=(a_1)_k\dots(a_m)_k
\end{equation*}
where $k$ is an integer or infinity. Further, we utilize
\begin{equation}\label{defhyp}
{}_s\phi_{s-1}\!\left[\begin{matrix}a_1,a_2,\dots,a_s\\
b_1,b_2,\dots,b_{s-1}\end{matrix}\,;q,z\right]:=
\sum _{k=0} ^{\infty}\frac {(a_1,a_2,\dots,a_s)_k}
{(q,b_1,\dots,b_{s-1})_k}z^k,
\end{equation}
and
\begin{equation}\label{defhypb}
{}_s\psi_s\!\left[\begin{matrix}a_1,a_2,\dots,a_s\\
b_1,b_2,\dots,b_s\end{matrix}\,;q,z\right]:=
\sum _{k=-\infty} ^{\infty}\frac {(a_1,a_2,\dots,a_s)_k}
{(b_1,b_2,\dots,b_s)_k}z^k,
\end{equation}
to denote the {\em basic hypergeometric ${}_s\phi_{s-1}$ series},
and the {\em bilateral basic hypergeometric ${}_s\psi_s$ series},
respectively. In \eqref{defhyp} or \eqref{defhypb}, $a_1,\dots,a_s$ are
called the {\em upper parameters}, $b_1,\dots,b_s$ the
{\em lower parameters}, $z$ is the {\em argument}, and 
$q$ the {\em base} of the series.
See \cite[p.~5 and p.~137]{GR} for the criteria
of when these series terminate, or, if not, when they converge. 
 
The classical theory of basic hypergeometric series contains
numerous summation and transformation formulae
involving ${}_s\phi_{s-1}$ or ${}_s\psi_s$ series.
Many of these summation theorems require
that the parameters satisfy the condition of being
either balanced and/or very-well-poised.
An ${}_s\phi_{s-1}$ basic hypergeometric series is called
{\em balanced} if $b_1\cdots b_{s-1}=a_1\cdots a_sq$ and $z=q$.
An ${}_s\phi_{s-1}$ series is {\em well-poised} if
$a_1q=a_2b_1=\cdots=a_sb_{s-1}$. An ${}_s\phi_{s-1}$ basic
hypergeometric series is called {\em very-well-poised}
if it is well-poised and if $a_2=-a_3=q\sqrt{a_1}$.
Note that the factor
\begin{equation*}
\frac {1-a_1q^{2k}}{1-a_1}
\end{equation*}
appears in a very-well-poised series.
The parameter $a_1$ is usually referred to as the
{\em special parameter} of such a series.
Similarly, a bilateral ${}_s\psi_s$ basic hypergeometric series is
well-poised if $a_1b_1=a_2b_2\cdots=a_sb_s$ and very-well-poised if,
in addition, $a_1=-a_2=qb_1=-qb_2$.

A standard reference for basic hypergeometric series
is Gasper and Rahman's texts~\cite{GR}.
In our computations in the subsequent sections
we frequently use some elementary identities of
$q$-shifted factorials, listed in \cite[Appendix~I]{GR}.

One of the most important theorems in the theory of basic
hypergeometric series is 
Jackson's~\cite{J} terminating very-well-poised balanced
${}_8\phi_7$ summation (cf.\ \cite[Eq.~(2.6.2)]{GR}):
\begin{multline}\label{87gl}
{}_8\phi_7\!\left[\begin{matrix}a,\,q\sqrt{a},-q\sqrt{a},b,c,d,
a^2q^{1+n}/bcd,q^{-n}\\
\sqrt{a},-\sqrt{a},aq/b,aq/c,aq/d,bcdq^{-n}/a,aq^{1+n}\end{matrix}\,;q,
q\right]\\
=\frac {(aq,aq/bc,aq/bd,aq/cd)_n}
{(aq/b,aq/c,aq/d,aq/bcd)_n}.
\end{multline}
This identity stands on the top of the classical hierarchy of summations
for basic hypergeometric series. Special cases include the terminating
and nonterminating very-well-poised ${}_6\phi_5$ summations,
the $q$-Pfaff--Saalsch\"utz summation, the $q$-Gau{\ss} summation,
the $q$-Chu--Vandermonde summation and the termininating and
nonterminating $q$-binomial theorem, see \cite{GR}.

\subsection{Multidimensional series}

$A_r$ (or, equivalently, $U(r+1)$) hypergeometric series were motivated
by the work of Biedenharn, Holman, and Louck~\cite{HBL}
in theoretical physics.
The theory of $A_r$ basic hypergeometric series
(or ``multiple basic hypergeometric series associated with the
root system $A_r$'', or ``associated with the unitary group $U(r+1)$''),
analogous to the classical theory of one-dimensional series,
has been developed originally by R.\ A.\ Gustafson, S.\ C.\ Milne,
and their co-workers, and later others (see \cite{BM,G1,G2,M1,M2,ML,Ro}
for a very small selection of papers in this area, neglecting a bunch
of other important references which already have grown vast in number).
Notably, several higher-dimensional extensions have been derived
(in each case) for the $q$-binomial theorem, $q$-Chu--Vandermonde
summation, $q$-Pfaff--Saalsch\"utz summation, Jackson's ${}_8\phi_7$
summation, Bailey's ${}_{10}\phi_9$ transformation, and
other important summation and transformation theorems.
See \cite{M3} for a survey on some of the main results and
techniques from the theory of $A_r$ basic hypergeometric series.

A recent major advance was the development of the theory of elliptic
hypergeometric series initiated by Frenkel and Turaev~\cite{FT}.
This ultimately lead to the study of elliptic hypergeometric series
associated with the root system $A_r$ (and other root systems),
see \cite{Ro}. Some of these developments are described in
Chapter~11 of Gasper and Rahman's texts~\cite{GR}.

A characteristic feature of $A_r$ series is that they
contain (the $A_r$-type product)
\begin{equation*}
\prod_{1\le i<j\le r}\frac{(1-q^{k_i-k_j}x_i/x_j)}{(1-x_i/x_j)}
\end{equation*}
as a factor in the summand (while they
should not, at the same time, contain factors that are characteristic
for other types of multiple series, such as for $C_r$ series,
see \cite{G2}). This characteristic feature is indeed shared by all the
multiple series considered in this paper which we therefore choose to
label as $A_r$ series.

When dealing with multivariable series,
we shall use the compact notations
\begin{equation*}
|\bk|:=k_1+\dots+k_r,\quad\text{where}\quad
\bk=(k_1,\dots,k_r),
\end{equation*}
and
\begin{equation*}
C:=c_1\cdots c_r,\qquad E:=e_1\cdots e_r.
\end{equation*}

We will need the following fundamental summation theorem by
Milne~\cite{M1}, originally obtained by
specializing an $A_r$ $q$-Whipple transformation
derived by partial fraction decompositions and functional
equations. For a simpler, more direct proof (of the elliptic
extension of Proposition~\ref{ar87}) see Rosengren~\cite{Ro},
which uses partial fraction decompositions and induction.

\begin{Proposition}[(Milne) An $A_r$ terminating very-well-poised
balanced $_8\phi_7$ summation]\label{ar87}
Let $a$, $b$, $c$, $d$ and $x_1,\dots,x_r$ be indeterminate,
let $n_1,\dots,n_r$ be nonnegative integers,
let $r\ge 1$, and suppose that
none of the denominators in \eqref{ar87gl} vanish. Then
\begin{multline}\label{ar87gl}
\underset{i=1,\dots,r}{\sum_{0\le k_i\le n_i}}
\prod_{i=1}^r\frac{(1-ax_iq^{k_i+|\bk|})}{(1-ax_i)}
\prod_{1\le i<j\le r}\frac{(1-q^{k_i-k_j}x_i/x_j)}{(1-x_i/x_j)}
\prod_{i,j=1}^r\frac{(q^{-n_j}x_i/x_j)_{k_i}}{(qx_i/x_j)_{k_i}}\\\times
\prod_{i=1}^r\frac{(ax_i)_{|\bk|}\,(dx_i,a^2x_iq^{1+|\bn|}/bcd)_{k_i}}
{(ax_iq^{1+|\bn|})_{|\bk|}\,(ax_iq/b,ax_iq/c)_{k_i}}\cdot
\frac{(b,c)_{|\bk|}}{(aq/d,bcdq^{-|\bn|}/a)_{|\bk|}}\,
q^{\sum_{i=1}^rik_i}\\=
\frac{(aq/bd,aq/cd)_{|\bn|}}{(aq/d,aq/bcd)_{|\bn|}}
\prod_{i=1}^r\frac{(ax_iq,ax_iq/bc)_{n_i}}{(ax_iq/b,ax_iq/c)_{n_i}}.
\end{multline}
\end{Proposition}

In Section~\ref{secn87} we apply inverse relations to
\eqref{ar87gl} to obtain a multivariable $_8\phi_7$ summation
of a different type by which we extend a theorem of
Bhatnagar~\cite[Thm.~3.6]{Bh}, whose result is stated as follows:

\begin{Proposition}[(Bhatnagar) An $A_r$ terminating
very-well-poised $_6\phi_5$ summation]\label{ar65}
Let $a$, $b$, $c$ and $x_1,\dots,x_r$ be indeterminate,
let $n_1,\dots,n_r$ be nonnegative integers,
let $r\ge 1$, and suppose that
none of the denominators in \eqref{ar65gl} vanish. Then
\begin{multline}\label{ar65gl}
\underset{i=1,\dots,r}{\sum_{0\le k_i\le n_i}}
\prod_{1\le i<j\le r}\frac{(1-q^{k_i-k_j}x_i/x_j)}{(1-x_i/x_j)}
\prod_{i,j=1}^r\frac{(q^{-n_j}x_i/x_j)_{k_i}}{(qx_i/x_j)_{k_i}}
\prod_{i=1}^r\frac{(c/x_i)_{|\bk|}\;x_i^{k_i}}
{(ax_iq/c)_{k_i}\,(c/x_i)_{|\bk|-k_i}}\\\times
\frac{(1-aq^{2|\bk|})}{(1-a)}
\frac{(a,b)_{|\bk|}}{(aq^{1+|\bn|},aq/b)_{|\bk|}}
\left(\frac{aq^{1+|\bn|}}{bc}\right)^{|\bk|}
q^{-e_2(\bk)+\sum_{i=1}^r(i-1)k_i}\\=
\frac{(aq)_{|\bn|}}{(aq/b)_{|\bn|}}
\prod_{i=1}^r\frac{(ax_iq/bc)_{n_i}}{(ax_iq/c)_{n_i}},
\end{multline}
where $e_2(\bk)$ is the second elementary symmetric function
of $\bk$.
\end{Proposition}

\begin{Remark}\label{rem}
The identity in Proposition~\ref{ar65} was derived in \cite{Bh}
by inverting an $A_r$ terminating balanced $_3\phi_2$ summation
from \cite{M2}. Bhatnagar calls his result a $U(r+1)$ summation
theorem (which in our terminology corresponds to $A_r$), since
the series in \eqref{ar65gl} looks very much like the usual
$U(r+1)$ (or $A_r$) series. We believe this classification to be
slightly inaccurate (but do not change it).
The series in \eqref{ar65gl} is apparently of a different type
than $A_r$ which (after $r\mapsto r+1$) can often be written
in form of an $\tilde{A}_r$ (or $SU(r)$) series.
It actually may be more accurate to call the series in \eqref{ar65gl}
(and the others in this paper, apart from the one appearing in
\eqref{ar87gl}) to be of type $A_1\times A_{r-1}$ (or similar),
due to the combination of two evidently different types of factors.
However, in lack of a more solid (say, representation-theoretic)
explanation, we (at least for now) refrain from labelling these series
as $A_1\times A_{r-1}$ and simply keep referring to them as $A_r$ series.
\end{Remark}

\begin{Remark}
We note that for $r\ge 2$ it is not possible to extend the
terminating summation in Proposition~\ref{ar65} by analytic
continuation to a nonterminating identity. This is due to the
appearance of the factor $q^{-e_2(\bk)}$ in the summand.
\end{Remark}

\section{Multidimensional matrix inversions}\label{secmmi}

Let $\mathbb Z$ denote the set of integers.
In the following,
we consider infinite lower-triangular $r$-dimensional matrices
$F=(f_{\bn\bk})_{\bn,\bk\in\mathbb Z^r}$ and
$G=(g_{\bn\bk})_{\bn,\bk\in\mathbb Z^r}$
(i.e., $f_{\bn\bk}=0$ unless $\bn\ge\bk$,
by which we mean $n_i\ge k_i$ for all $i=1,\dots,r$), and
infinite sequences $(a_\bn)_{\bn\in\mathbb Z^r}$ and
$(b_\bn)_{\bn\in\mathbb Z^r}$.

The matrix $F$ is said to be the {\em inverse} of $G$,
if and only if the following orthogonality
relation holds:
\begin{equation}\label{orthrel}
\sum _{\bn\ge\bk\ge\bl}
f_{\bn\bk}g_{\bk\bl}=
\delta_{\bn\bl}\qquad\qquad
\text {for all}\quad \bn,\bl\in\mathbb Z^r,
\end{equation}
where $\delta_{\bn\bl}$ is the usual Kronecker delta.
Since $F$ and $G$ are lower-triangular, the sum in \eqref{orthrel}
is finite and also the dual relation, with the roles of $F$ and $G$
being interchanged, must hold at the same time.

It follows readily from the orthogonality relation \eqref{orthrel}
that
\begin{subequations}\label{invrel}
\begin{equation}\label{invrel1}
\sum _{{\mathbf 0}\le\bk\le\bn}
f_{\bn\bk}a_\bk=b_\bn
\qquad\qquad\text{for all $\bn\in\mathbb Z^r$,}
\end{equation}
if and only if
\begin{equation}\label{invrel2}
\sum _{{\mathbf 0}\le\bl\le\bk}
g_{\bk\bl}b_\bl=a_\bk
\qquad\qquad\text{for all $\bk\in\mathbb Z^r$.}
\end{equation}
\end{subequations}

Inverse relations are a powerful tool for proving or deriving
identities. For instance, given an identity in the form \eqref{invrel1},
we can immediately deduce \eqref{invrel2}, which may possibly be a new
identity. It is exactly this variant of multiple inverse relations
which we apply in the derivation of Theorem~\ref{n87}.

One of the main results of \cite{S1} (see Thm.~3.1 therein)
was the following explicit multidimensional matrix inverse
which reduces to Krattenthaler's matrix inverse~\cite{Kr} for $r=1$.

\begin{Proposition}[A general $A_r$ matrix inverse]\label{mmi}
Let $(a_t)_{t\in\mathbb Z}$ and $(c_j(t))_{t\in\mathbb
Z}$, $1\leq j\leq r$, be arbitrary sequences of scalars. Then the
lower-triangular $r$-dimensional matrices 
\begin{subequations}\label{mirs}
\begin{equation}
f_{\bf nk}=\frac{\prod_{t=|{\bf k}|}^{|{\bf n}|-1}
\left((1-a_tc_1(k_1)\dotsm c_r(k_r))
\prod_{j=1}^r(a_t-c_j(k_j))\right)} 
{\prod_{i=1}^r\prod_{t=k_i+1}^{n_i}\left((1-c_i(t)c_1(k_1)\dotsm
c_r(k_r))\prod_{j=1}^r(c_i(t)-c_j(k_j))\right)},
\end{equation}
and
\begin{equation}
\begin{split}g_{\bf kl}&=\prod_{1\leq i<j\leq r}
\frac{(c_i(l_i)-c_j(l_j))}
{(c_i(k_i)-c_j(k_j))}\,\cdot
\frac{(1-a_{|{\bf l}|}c_1(l_1)\dotsm c_r(l_r))}
{(1-a_{|{\bf k}|}c_1(k_1)\dotsm c_r(k_r))}
\prod_{j=1}^r\frac{(a_{|{\bf l}|}-c_j(l_j))}
{(a_{|{\bf k}|}-c_j(k_j))}\\
&\quad\times\frac{\prod_{t=|{\bf l}|+1}^{|{\bf k}|}
\left((1-a_tc_1(k_1)\dotsm c_r(k_r))
\prod_{j=1}^r(a_t-c_j(k_j))\right)} 
{\prod_{i=1}^r\prod_{t=l_i}^{k_i-1}\left((1-c_i(t)c_1(k_1)\dotsm c_r(k_r))
\prod_{j=1}^r(c_i(t)-c_j(k_j))\right)}
\end{split}
\end{equation}
\end{subequations}
are mutually inverse.
\end{Proposition}

For an elliptic extension of the above, see \cite{RS}.
Some important special cases of Proposition~\ref{mmi} include
Bhatnagar and Milne's $A_r$ matrix inverse \cite[Thm.~3.48]{BM}
(take $c_i(t)\mapsto x_iq^{t}$, for $i=1,\dots,r$),
the author's $D_r$ matrix inverse \cite[Thm.~5.11]{S1}
(take $a_t=0$ and $c_i(t)\mapsto x_iq^{t}+q^{-t}/x_i$, for $i=1,\dots,r$),
and some other (non-hypergeometric) matrix inverses
considered in \cite[App.~A]{S2}.

The following corollary of Proposition~\ref{mmi}
has not yet been explicitly stated, nor applied.
(This is maybe surprising as it also can be derived
from Bhatnagar and Milne's matrix inverse.)
It is readily obtained from \eqref{mirs} by letting
$a_t\mapsto (b/x_1\dots x_r)^{\frac 1{r+1}}aq^{t}$ and
$c_i(t)\mapsto (b/x_1\dots x_r)^{\frac 1{r+1}}x_iq^{t}$,
for $i=1,\dots,r$, followed by some simplifications including
\begin{multline}
\prod_{1\le i<j\le r}\frac{(1-q^{k_i-k_j}x_i/x_j)}
{(1-q^{l_i-l_j}x_i/x_j)}
\prod_{i,j=1}^r\frac{(q^{l_i-k_j}x_i/x_j)_{k_i-l_i}}
{(q^{1+l_i-l_j}x_i/x_j)_{k_i-l_i}}\\
=(-1)^{|{\mathbf k}|-|{\mathbf l}|}\,
q^{-\binom{|{\mathbf k}|-|{\mathbf l}|}2-\sum_{i=1}^{r}i(k_i-l_i)},
\end{multline}
the latter of which is equivalent to Lemma~4.3 of \cite{M2}
and is typical for dealing with $A_r$ series.

\begin{Corollary}[An $A_r$ matrix inverse]\label{mmic}
Let $a$, $b$ and $x_1\dots,x_r$ be indeterminates.
Then the lower-triangu\-lar $r$-dimensional matrices 
\begin{subequations}\label{migl}
\begin{equation}\label{migla}
f_{\bn\bk}=
\frac{(abq^{2|\bk|})_{|\bn|-|\bk|}
\prod_{i=1}^r(aq^{|\bk|-k_i}/x_i)_{|\bn|-|\bk|}}
{\prod_{i=1}^r(bx_iq^{1+k_i+|\bk|})_{n_i-k_i}
\prod_{i,j=1}^r(q^{1+k_i-k_j}x_i/x_j)_{n_i-k_i}},
\end{equation}
\begin{align}\label{miglb}\notag
g_{\bk\bl}={}&
(-1)^{|\bk|-|\bl|}q^{\binom{|\bk|-|\bl|}2}
\frac{(1-abq^{2|\bl|})}{(1-abq^{2|\bk|})}
\prod_{i=1}^r\frac{(1-aq^{|\bl|-l_i}/x_i)}{(1-aq^{|\bk|-k_i}/x_i)}\\
&\times
\frac{(abq^{1+|\bl|+|\bk|})_{|\bk|-|\bl|}
\prod_{i=1}^r(aq^{1+|\bl|-k_i}/x_i)_{|\bk|-|\bl|}}
{\prod_{i=1}^r(bx_iq^{l_i+|\bk|})_{k_i-l_i}
\prod_{i,j=1}^r(q^{1+l_i-l_j}x_i/x_j)_{k_i-l_i}}.
\end{align}
\end{subequations}
are mutually inverse.
\end{Corollary}

Corollary~\ref{mmic} constitutes a multivariable extension
of Bressoud's matrix inverse~\cite{Br} which he extracted
directly from the terminating very-well-poised ${}_6\phi_5$
summation. Bressoud's matrix inverse underlies the
WP-Bailey lemma \cite{A2} which is a generalization
of the classical Bailey lemma \cite{A1}, both powerful tools
for deriving (chains of) identities.
Other multivariable extensions of Bressoud's matrix inverse
have been derived (for type $A$) by Milne~\cite[Thm.~3.41]{M2},
(for type $C$) by Lilly and Milne~\cite[2nd Remark after Thm.~2.11]{LM},
(for type $D$) by the author~\cite[Thm.~5.11]{S1},
(of Carlitz type, `twisted') by Krattenthaler and the
author~\cite[Eqs.~(6.4)/(6.5)]{KS},
(related to $A_{n-1}$ Macdonald polynomials) by Lassalle and the
author~\cite[Thm.~2.7]{LS},
and (related to an elliptic extension of
$BC_n$ Koornwinder--Macdonald polynomials) by Rains~\cite[Cor.\ 4.3]{Ra}
and by Coskun and Gustafson~\cite[Eq.~(4.16)]{CG},
and possibly others (of which the current author is not aware of).

\section{New $A_r$ terminating very-well-poised
${}_8\phi_7$ summations}\label{secn87}

We now combine the multidimensional matrix inverse in
Corollary~\ref{mmic} with Milne's $A_r$ ${}_8\phi_7$ summation
in Proposition~\ref{ar87} to deduce a new multivariable ${}_8\phi_7$
summation theorem.

In particular, we have \eqref{invrel1} by the
$(a,b,c)\mapsto (b,c,abq^{|\bn|})$ case of Proposition~\ref{ar87},
where
\begin{align*}\notag
a_{\bk}={}&\frac{(ab)_{2|\bk|}\,(c)_{|\bk|}}{(acd,bq/d)_{|\bk|}}
\prod_{i=1}^r\frac{(bx_i)_{|\bk|}\,(a/x_i)_{|\bk|-k_i}\,
(dx_i,bx_iq/acd)_{k_i}}{(bx_i)_{k_i+|\bk|}\,(bx_iq/c)_{k_i}}\\&\times
q^{\binom{|\bk|}2-\sum_{i=1}^r\binom{k_i}2}a^{|\bk|}
\prod_{i=1}^rx_i^{-k_i}\prod_{i,j=1}^r(qx_i/x_j)_{k_i}^{-1}
\end{align*}
and
\begin{equation*}
b_{\bn}=\frac{(ab,ad,bq/cd)_{|\bn|}}{(acd,bq/d)_{|\bn|}}
\prod_{i=1}^r\frac{(ac/x_i)_{|\bn|}\,(a/x_i)_{|\bn|-n_i}}
{(bx_iq/c)_{n_i}\,(ac/x_i)_{|\bn|-n_i}}
\prod_{i,j=1}^r(qx_i/x_j)_{n_i}^{-1},
\end{equation*}
and $f_{\bn\bk}$ as in \eqref{migla}. Therefore we must have
\eqref{invrel2} with the above sequences $b_{\bl}$, $a_{\bk}$,
and $g_{\bk\bl}$ as in \eqref{miglb}. In explicit terms this
gives (after simplifications and the substitutions
$(a,c,d,x_i,k_i,l_i)\mapsto(a/b,aq/bc,b^2/a,a^2qx_i/b^2cd,n_i,k_i)$,
$i=1,\dots,r$) the following new multivariable extension of \eqref{87gl}:

\begin{Theorem}[An $A_r$ terminating very-well-poised
balanced $_8\phi_7$ summation]\label{n87}
Let $a$, $b$, $c$, $d$ and $x_1,\dots,x_r$ be indeterminate,
let $n_1,\dots,n_r$ be nonnegative integers, let $r\ge 1$, and suppose that
none of the denominators in \eqref{n87gl} vanish. Then
\begin{multline}\label{n87gl}
\underset{i=1,\dots,r}{\sum_{0\le k_i\le n_i}}
\prod_{1\le i<j\le r}\frac{(1-q^{k_i-k_j}x_i/x_j)}{(1-x_i/x_j)}
\prod_{i,j=1}^r\frac{(q^{-n_j}x_i/x_j)_{k_i}}{(qx_i/x_j)_{k_i}}\\\times
\prod_{i=1}^r\frac{(bcd/ax_i)_{|\bk|-k_i}\,(d/x_i)_{|\bk|}\,
(a^2x_iq^{1+|\bn|}/bcd)_{k_i}}
{(d/x_i)_{|\bk|-k_i}\,(bcdq^{-n_i}/ax_i)_{|\bk|}\,(ax_iq/d)_{k_i}}\\\times
\frac{(1-aq^{2|\bk|})}{(1-a)}
\frac{(a,b,c)_{|\bk|}}{(aq^{1+|\bn|},aq/b,aq/c)_{|\bk|}}\,
q^{\sum_{i=1}^rik_i}\\=
\frac{(aq,aq/bc)_{|\bn|}}{(aq/b,aq/c)_{|\bn|}}
\prod_{i=1}^r\frac{(ax_iq/bd,ax_iq/cd)_{n_i}}{(ax_iq/d,ax_iq/bcd)_{n_i}}.
\end{multline}
\end{Theorem}

By a polynomial argument, this is equivalent to the following result.

\begin{Corollary}[An $A_r$ terminating very-well-poised
balanced $_8\phi_7$ summation]\label{n87c}
Let $a$, $b$, $c_1,\dots,c_r$, $d$ and $x_1,\dots,x_r$ be indeterminate,
let $N$ be a nonnegative integer, let $r\ge 1$, and suppose that
none of the denominators in \eqref{n87cgl} vanish. Then
\begin{multline}\label{n87cgl}
\underset{0\le|\bk|\le N}{\sum_{k_1,\dots,k_r\ge 0}}
\prod_{1\le i<j\le r}\frac{(1-q^{k_i-k_j}x_i/x_j)}{(1-x_i/x_j)}
\prod_{i,j=1}^r\frac{(c_jx_i/x_j)_{k_i}}{(qx_i/x_j)_{k_i}}\\\times
\prod_{i=1}^r\frac{(bdq^{-N}/ax_i)_{|\bk|-k_i}\,(d/x_i)_{|\bk|}\,
(a^2x_iq^{1+N}/bCd)_{k_i}}
{(d/x_i)_{|\bk|-k_i}\,(bc_idq^{-N}/ax_i)_{|\bk|}\,(ax_iq/d)_{k_i}}\\\times
\frac{(1-aq^{2|\bk|})}{(1-a)}
\frac{(a,b,q^{-N})_{|\bk|}}{(aq/C,aq/b,aq^{1+N})_{|\bk|}}\,
q^{\sum_{i=1}^rik_i}\\=
\frac{(aq,aq/bC)_{N}}{(aq/b,aq/C)_{N}}
\prod_{i=1}^r\frac{(ax_iq/bd,ax_iq/c_id)_{N}}{(ax_iq/d,ax_iq/bc_id)_{N}},
\end{multline}
where $C=c_1\cdots c_r$.
\end{Corollary}

\begin{proof}
First we write the right side of \eqref{n87cgl} as quotient of infinite
products using \eqref{ipr}. Then by the $c=q^{-N}$ case of
Theorem~\ref{n87} it follows that the identity \eqref{n87cgl}
holds for $c_j=q^{-n_j}$, $j=1,\dots,r$.
By clearing out denominators in \eqref{n87cgl},
we get a polynomial equation in $c_1$, which is true for $q^{-n_1}$,
$n_1=0,1,\dots$. Thus we obtain an identity in $c_1$.
By carrying out this process for $c_2,c_3,\dots,c_r$ also,
we obtain Corollary~\ref{n87c}.
\end{proof}

\begin{Remark}
For $b\to\infty$, Theorem~\ref{n87} and Corollary~\ref{n87c} reduce
(after relabeling of parameters) to the two respective terminating
${}_6\phi_5$ summations in Bhatnagar~\cite[Thms.~3.6 and 3.7]{Bh},
of which the first is displayed in Proposition~\ref{ar65}.
\end{Remark}

The (equivalent) $b\mapsto a^2q^{1+N}/bCd$ case of \eqref{n87cgl}
appears to be particularly useful:
\begin{multline}\label{n87ogl}
\underset{0\le|\bk|\le N}{\sum_{k_1,\dots,k_r\ge 0}}
\prod_{1\le i<j\le r}\frac{(1-q^{k_i-k_j}x_i/x_j)}{(1-x_i/x_j)}
\prod_{i,j=1}^r\frac{(c_jx_i/x_j)_{k_i}}{(qx_i/x_j)_{k_i}}\\\times
\prod_{i=1}^r\frac{(aq/bCx_i)_{|\bk|-k_i}\,(d/x_i)_{|\bk|}\,
(bx_i)_{k_i}}
{(d/x_i)_{|\bk|-k_i}\,(ac_iq/bCx_i)_{|\bk|}\,(ax_iq/d)_{k_i}}\\\times
\frac{(1-aq^{2|\bk|})}{(1-a)}
\frac{(a,a^2q^{1+N}/bCd,q^{-N})_{|\bk|}}
{(aq/C,bCdq^{-N}/a,aq^{1+N})_{|\bk|}}\,
q^{\sum_{i=1}^rik_i}\\=
\frac{(aq,aq/bd)_{N}}{(aq/C,aq/bCd)_{N}}
\prod_{i=1}^r\frac{(aq/bCx_i,ax_iq/c_id)_{N}}{(ax_iq/d,ac_iq/bCx_i)_{N}},
\end{multline}
where $C=c_1\cdots c_r$.

Note that in the summand of the series on the left-hand side of
\eqref{n87ogl} the terminating integer $N$ appears only within
factors depending on $|\bk|$. This makes it particularly convenient
to combine \eqref{n87ogl} with sums depending $N$ to obtain further
results such as a new multivariable ${}_{10}\phi_9$ transformation.
(See \cite{RS} for details, in the more general setting of elliptic
hypergeometric series.) 

\section{New $A_r$ terminating and nonterminating
very-well-poised ${}_6\phi_5$ summations}\label{secn65}

In \eqref{n87ogl} we now let $N\to\infty$ (while appealing to
Tannery's theorem (cf.\ ~\cite{Bw}) for justification of taking
term-wise limits) and obtain the following result:

\begin{Corollary}[An $A_r$ nonterminating very-well-poised
${}_6\phi_5$ summation]\label{n65}
Let $a$, $b$, $c_1,\dots,c_r$, $d$ and $x_1,\dots,x_r$ be indeterminate,
let $r\ge 1$, and suppose that
none of the denominators in \eqref{n65gl} vanish. Then
\begin{multline}\label{n65gl}
\sum_{k_1,\dots,k_r\ge 0}
\prod_{1\le i<j\le r}\frac{(1-q^{k_i-k_j}x_i/x_j)}{(1-x_i/x_j)}
\prod_{i,j=1}^r\frac{(c_jx_i/x_j)_{k_i}}{(qx_i/x_j)_{k_i}}
\\\times
\prod_{i=1}^r\frac{(aq/bCx_i)_{|\bk|-k_i}\,(d/x_i)_{|\bk|}\,
(bx_i)_{k_i}}
{(d/x_i)_{|\bk|-k_i}\,(ac_iq/bCx_i)_{|\bk|}\,(ax_iq/d)_{k_i}}\\\times
\frac{(1-aq^{2|\bk|})}{(1-a)}\frac{(a)_{|\bk|}}{(aq/C)_{|\bk|}}
\left(\frac{aq}{bCd}\right)^{|\bk|}
q^{\sum_{i=1}^r(i-1)k_i}\\=
\frac{(aq,aq/bd)_\infty}{(aq/C,aq/bCd)_\infty}
\prod_{i=1}^r\frac{(aq/bCx_i,ax_iq/c_id)_\infty}
{(ax_iq/d,ac_iq/bCx_i)_\infty},
\end{multline}
provided $|aq/bCd|<1$, where $C=c_1\cdots c_r$.
\end{Corollary}

An immediate consequence of Corollary~\ref{n65}
(obtained by letting $c_i=q^{-n_i}$, $i=1,\dots,r$, and
$d\mapsto c$, in \eqref{n65gl}) is the following terminating summation:

\begin{Corollary}[An $A_r$ terminating very-well-poised
${}_6\phi_5$ summation]\label{n65t}
Let $a$, $b$, $c$ and $x_1,\dots,x_r$ be indeterminate,
let $n_1,\dots,n_r$ be nonnegative integers,
let $r\ge 1$, and suppose that
none of the denominators in \eqref{n65tgl} vanish. Then
\begin{multline}\label{n65tgl}
\underset{i=1,\dots,r}{\sum_{0\le k_i\le n_i}}
\prod_{1\le i<j\le r}\frac{(1-q^{k_i-k_j}x_i/x_j)}{(1-x_i/x_j)}
\prod_{i,j=1}^r\frac{(q^{-n_j}x_i/x_j)_{k_i}}{(qx_i/x_j)_{k_i}}\\\times
\prod_{i=1}^r\frac{(aq^{1+|\bn|}/bx_i)_{|\bk|-k_i}\,(c/x_i)_{|\bk|}\,
(bx_i)_{k_i}}
{(c/x_i)_{|\bk|-k_i}\,(aq^{1+|\bn|-n_i}/bx_i)_{|\bk|}\,
(ax_iq/c)_{k_i}}\\\times
\frac{(1-aq^{2|\bk|})}{(1-a)}\frac{(a)_{|\bk|}}{(aq^{1+|\bn|})_{|\bk|}}
\left(\frac{aq^{1+|\bn|}}{bc}\right)^{|\bk|}
q^{\sum_{i=1}^r(i-1)k_i}\\=
(aq,aq/bc)_{|\bn|}
\prod_{i=1}^r\frac{(aq/bx_i)_{|\bn|-n_i}}{(aq/bx_i)_{|\bn|}
(ax_iq/c)_{n_i}}.
\end{multline}
\end{Corollary}

Corollary~\ref{n65t} can also be obtained from
Theorem~\ref{n87} by first replacing $b$ by $a^2q^{1+|\bn|}/bcd$,
then letting $c\to\infty$, followed by relabeling $d\mapsto c$.

On the other hand, we can also deduce the nonterminating ${}_6\phi_5$
summation in Corollary~\ref{n65} from the terminating sum in
Corollary~\ref{n65t} by analytic continuation (by which one can
pretend to avoid the explicit application of Tannery's theorem
in the derivation of Corollary~\ref{n65}, however, implicitly such
an application is needed to show the analyticity of the series),
in the form of a repeated application of a variant of Ismail's
argument \cite{Is}.
Indeed, both sides of the multiple series identity in
\eqref{n65gl} are analytic in each of the parameters
$1/c_1,\dots,1/c_r$ in a domain around the origin
(see the following paragraph).
Now, the identity is true for $1/c_1=q^{n_1},
1/c_2=q^{n_2},\dots,$ and $1/c_r=q^{n_r}$, by the
$A_r$ terminating ${}_6\phi_5$ summation in
Corollary~\ref{n65t}. This holds for
all $n_1,\dots,n_r\ge 0$.
Since $\lim_{n_1\to\infty}q^{n_1}=0$ is an interior point
in the domain of analyticity of $1/c_1$, by the identity theorem
we obtain an identity for general $1/c_1$.
By iterating this argument for $1/c_2,\dots,1/c_r$, we establish
\eqref{n65gl} for general $1/c_1,\dots,1/c_r$
(i.e., for general $c_1,\dots,c_r$).

What remains to be shown is the claim that both sides of
\eqref{n65gl} are analytic in the parameters $1/c_1,\dots,1/c_r$.
This is accomplished by multiple applications
of the $q$-binomial theorem (cf.\ \cite{GR})
\begin{equation}\label{qbin}
\sum_{l\ge 0}\frac{(a)_l}{(q)_l}z^j=\frac{(az)_\infty}{(z)_\infty},
\qquad\qquad\text{where}\quad |z|<1,
\end{equation}
to expand both sides of the identity as a convergent
multiple power series in $1/c_1,\dots,1/c_r$.
We leave the details, similar in nature to those in the proof of
Theorem~\ref{n66}, to the reader.

\section{A new $A_r$ very-well-poised ${}_6\psi_6$ summation}\label{secn66}

Having Corollary~\ref{n65}, we are ready to prove the following
multivariable extension of Bailey's very-well-poised
${}_6\psi_6$ summation formula in \eqref{66gl}:
 
\begin{Theorem}[An $A_r$ very-well-poised ${}_6\psi_6$ summation]\label{n66}
Let $a$, $b$, $c_1,\dots,c_r$, $d$, $e_1,\dots,e_r$ and $x_1,\dots,x_r$
be indeterminate, let $r\ge 1$, and suppose that
none of the denominators in \eqref{n66gl} vanish. Then
\begin{multline}\label{n66gl}
\sum_{-\infty\le k_1,\dots,k_r\le\infty}
\prod_{1\le i<j\le r}\frac{(1-q^{k_i-k_j}x_i/x_j)}{(1-x_i/x_j)}
\prod_{i,j=1}^r\frac{(c_jx_i/x_j)_{k_i}}{(ax_iq/e_jx_j)_{k_i}}
\\\times
\prod_{i=1}^r\frac{(aq/bCx_i)_{|\bk|-k_i}\,(dE/a^{r-1}e_ix_i)_{|\bk|}\,
(bx_i)_{k_i}}
{(dE/a^rx_i)_{|\bk|-k_i}\,(ac_iq/bCx_i)_{|\bk|}\,(ax_iq/d)_{k_i}}\\\times
\frac{(1-aq^{2|\bk|})}{(1-a)}\frac{(E/a^{r-1})_{|\bk|}}{(aq/C)_{|\bk|}}
\left(\frac{a^{r+1}q}{bCdE}\right)^{|\bk|}
q^{\sum_{i=1}^r(i-1)k_i}\\=
\frac{(aq,q/a,aq/bd)_\infty}{(aq/C,a^{r+1}q/bCdE,a^{r-1}q/E)_\infty}
\prod_{i,j=1}^r\frac{(qx_i/x_j,ax_iq/c_ie_jx_j)_\infty}
{(qx_i/c_ix_j,ax_iq/e_jx_j)_\infty}\\\times
\prod_{i=1}^r\frac{(a^rx_iq/dE,aq/be_ix_i,aq/bCx_i,ax_iq/c_id)_\infty}
{(a^{r-1}e_ix_iq/dE,q/bx_i,ax_iq/d,ac_iq/bCx_i)_\infty},
\end{multline}
provided $|aq^{r+1}/bCdE|<1$, where $C=c_1\cdots c_r$ and $E=e_1\cdots e_r$.
\end{Theorem}

Clearly, Theorem~\ref{n66} reduces to Corollary~\ref{n65} for
$e_1=e_2=\dots=e_r=a$. A different very-well-poised ${}_6\psi_6$
summation for the root system $A_r$ was given by Gustafson, see
\cite{G1} and \cite{G2}. Other multivariable very-well-poised
${}_6\psi_6$ summations are listed by Ito~\cite{I2}.

\begin{proof}[Proof of Theorem~\ref{n66}]
We apply Ismail's argument~\cite{Is} (see also \cite{AI})
successively to the parameters $1/e_1,\dots,1/e_r$ using the $A_r$
nonterminating ${}_6\phi_5$ summation in Corollary~\ref{n65}.
The multiple series identity on the left-hand side of \eqref{n66gl}
is analytic in each of the parameters
$1/e_1,\dots,1/e_r$ in a domain around the origin
(which is not difficult to verify, see further below).
Now, the identity is true for $1/e_1=q^{m_1}/a,
1/e_2=q^{m_2}/a,\dots,$ and $1/e_r=q^{m_r}/a$, by
Corollary~\ref{n65} (see the next paragraph for the details).
This holds for all $m_1,\dots,m_n\ge 0$.
Since $\lim_{m_1\to\infty}q^{m_1}/a=0$ is an interior point
in the domain of analyticity of $1/e_1$, by the identity theorem,
we obtain an identity for general $1/e_1$.
By iterating this argument for $1/e_2,\dots,1/e_r$, we establish
\eqref{n66gl} for general $1/e_1,\dots,1/e_r$
(i.e., for general $e_1,\dots,e_r$).

The details are displayed as follows. Setting $1/e_i=q^{m_i}/a$, for
$i=1,\dots,r$, the left-hand side of \eqref{n66gl} becomes
\begin{multline}\label{longl1}
\underset{i=1,\dots,r}{\sum_{-m_i\le k_i\le\infty}}
\prod_{1\le i<j\le r}\frac{(1-q^{k_i-k_j}x_i/x_j)}{(1-x_i/x_j)}
\prod_{i,j=1}^r\frac{(c_jx_i/x_j)_{k_i}}{(q^{1+m_j}x_i/x_j)_{k_i}}
\\\times
\prod_{i=1}^r\frac{(aq/bCx_i)_{|\bk|-k_i}\,(dq^{m_i-|\bm|}/x_i)_{|\bk|}\,
(bx_i)_{k_i}}
{(dq^{-|\bm|}/x_i)_{|\bk|-k_i}\,(ac_iq/bCx_i)_{|\bk|}\,(ax_iq/d)_{k_i}}\\\times
\frac{(1-aq^{2|\bk|})}{(1-a)}\frac{(aq^{-|\bm|})_{|\bk|}}{(aq/C)_{|\bk|}}
\left(\frac{aq^{1+|\bm|}}{bCd}\right)^{|\bk|}
q^{\sum_{i=1}^r(i-1)k_i}.
\end{multline}
We shift the summation indices in \eqref{longl1} by
$k_i\mapsto k_i-m_i$, for $i=1,\dots,r$ and obtain
\begin{multline*}
\prod_{1\le i<j\le r}\frac{(1-q^{m_j-m_i}x_i/x_j)}{(1-x_i/x_j)}
\prod_{i,j=1}^r\frac{(c_jx_i/x_j)_{-m_i}}{(q^{1+m_j}x_i/x_j)_{-m_i}}
\\\times
\prod_{i=1}^r\frac{(aq/bCx_i)_{m_i-|\bm|}\,
(dq^{m_i-|\bm|}/x_i)_{-|\bm|}\,(bx_i)_{-m_i}}
{(dq^{-|\bm|}/x_i)_{m_i-|\bm|}\,(ac_iq/bCx_i)_{-|\bm|}\,
(ax_iq/d)_{-m_i}}\\\times
\frac{(1-aq^{-2|\bm|})}{(1-a)}
\frac{(aq^{-|\bm|})_{-|\bm|}}{(aq/C)_{-|\bm|}}
\left(\frac{aq^{1+|\bm|}}{bCd}\right)^{-|\bm|}
q^{-\sum_{i=1}^r(i-1)m_i}\\\times
\sum_{k_1,\dots k_r\ge 0}
\prod_{1\le i<j\le r}\frac{(1-q^{m_j-m_i+k_i-k_j}x_i/x_j)}
{(1-q^{m_j-m_i}x_i/x_j)}
\prod_{i,j=1}^r\frac{(c_jq^{-m_i}x_i/x_j)_{k_i}}
{(q^{1+m_j-m_i}x_i/x_j)_{k_i}}\\\times
\prod_{i=1}^r\frac{(aq^{1+m_i-|\bm|}/bCx_i)_{|\bk|-k_i}\,
(dq^{m_i-2|\bm|}/x_i)_{|\bk|}\,(bx_iq^{-m_i})_{k_i}}
{(dq^{m_i-2|\bm|}/x_i)_{|\bk|-k_i}\,
(ac_iq^{1-|\bm|}/bCx_i)_{|\bk|}\,(ax_iq^{1-m_i}/d)_{k_i}}\\\times
\frac{(1-aq^{-2|\bm|+2|\bk|})}{(1-aq^{-2|\bm|})}
\frac{(aq^{-2|\bm|})_{|\bk|}}{(aq^{1-|\bm|}/C)_{|\bk|}}
\left(\frac{aq^{1+|\bm|}}{bCd}\right)^{|\bk|}
q^{\sum_{i=1}^r(i-1)k_i}.
\end{multline*}
Next, we apply the $a\mapsto aq^{-2|\bm|}$,
$x_i\mapsto x_iq^{-m_i}$, $c_i\mapsto c_iq^{-m_i}$, $i=1,\dots,r$,
$d\mapsto dq^{-2|\bm|}$ case of Corollary~\ref{n65}, and obtain
\begin{multline}\label{long}
\prod_{1\le i<j\le r}\frac{(1-q^{m_j-m_i}x_i/x_j)}{(1-x_i/x_j)}
\prod_{i,j=1}^r\frac{(c_jx_i/x_j)_{-m_i}}{(q^{1+m_j}x_i/x_j)_{-m_i}}
\\\times
\prod_{i=1}^r\frac{(aq/bCx_i)_{m_i-|\bm|}\,
(dq^{m_i-|\bm|}/x_i)_{-|\bm|}\,(bx_i)_{-m_i}}
{(dq^{-|\bm|}/x_i)_{m_i-|\bm|}\,(ac_iq/bCx_i)_{-|\bm|}\,
(ax_iq/d)_{-m_i}}\\\times
\frac{(1-aq^{-2|\bm|})}{(1-a)}
\frac{(aq^{-|\bm|})_{-|\bm|}}{(aq/C)_{-|\bm|}}
\left(\frac{aq^{1+|\bm|}}{bCd}\right)^{-|\bm|}
q^{-\sum_{i=1}^r(i-1)m_i}\\\times
\frac{(aq^{1-2|\bm|},aq/bd)_\infty}
{(aq^{1-|\bm|}/C,aq^{1+|\bm|}/bCd)_\infty}
\prod_{i=1}^r\frac{(aq^{1+m_i-|\bm|}/bCx_i,ax_iq/c_id)_\infty}
{(ax_iq^{1-m_i}/d,ac_iq^{1-|\bm|}/bCx_i)_\infty}.
\end{multline}
Finally, we apply several elementary identities from \cite[App.~I]{GR},
and the $n\mapsto r$, $y_i\mapsto-m_i$, $i=1,\dots,r$,
case of \cite[Lem.~3.12]{M2}, specifically
\begin{multline}
\prod_{i,j=1}^r(qx_i/x_j)_{m_j-m_i}=
(-1)^{(r-1)|{\bm}|}\,q^{-\binom{|\bm|+1}2+
r\sum_{i=1}^r\binom{m_i+1}2}\,q^{-\sum_{i=1}^r(i-1)m_i}\\\times
\prod_{i=1}^rx_i^{|{\mathbf m}|-rm_i}
\prod_{1\le i<j\le r}\frac{(1-q^{m_j-m_i}x_i/x_j)}{(1-x_i/x_j)},
\end{multline}
to transform the expression obtained in \eqref{long} to
\begin{multline*}
\frac{(aq,q/a,aq/bd)_\infty}{(aq/C,aq^{1+|\bm|}/bCd,q^{1+|\bm|}/a)_\infty}
\prod_{i,j=1}^r\frac{(qx_i/x_j,q^{1+m_j}x_i/c_ix_j)_\infty}
{(qx_i/c_ix_j,q^{1+m_j}x_i/x_j)_\infty}\\\times
\prod_{i=1}^r\frac{(x_iq^{1+|\bm|}/d,q^{1+m_i}/bx_i,
aq/bCx_i,ax_iq/c_id)_\infty}
{(x_iq^{1+|\bm|-m_i}/d,q/bx_i,ax_iq/d,ac_iq/bCx_i)_\infty},
\end{multline*}
which is exactly the $1/e_i=q^{m_i}/a$, $i=1,\dots,r$, case of the
right-hand side of \eqref{n66gl}.

We still need to show that both sides of \eqref{n66gl} are analytic
in $1/e_1,\dots,1/e_r$ around the origin, i.e., that both sides
can be expanded as convergent multiple powers of
the variables $1/e_1,\dots,1/e_r$. This is easily achieved
by multiple use of the $q$-binomial theorem \eqref{qbin}.
For the right-hand side the claim is immediate
(because all we need is to multiply expressions of the form
$(z)_\infty=\sum(-1)^lq^{\binom l2}z^l/(q)_l$
and $(z)_\infty^{-1}=\sum z^l/(q)_l$).
At the left-hand side we manipulate those factors in the summand of the
series which involve $1/e_1,\dots,1/e_r$ (in order to obtain a
convergent multiple power series).
In particular, we use
\begin{multline}\label{auxid}
\prod_{i,j=1}^r\frac 1{(ax_iq/e_jx_j)_{k_i}}
\prod_{i=1}^r\frac{(dE/a^{r-1}e_ix_i)_{|\bk|}}{(dE/a^rx_i)_{|\bk|-k_i}}
\cdot(E/a^{r-1})_{|\bk|}E^{-|\bk|}\\
=\prod_{i,j=1}^r\frac {(ax_iq^{1+k_i}/e_jx_j)_\infty}{(ax_iq/e_jx_j)_\infty}
\prod_{i=1}^r\frac{(a^{r-1}e_ix_iq^{1-|\bk|}/dE)_{|\bk|}}
{(a^rx_iq^{1+k_i-|\bk|}/dE)_{|\bk|-k_i}}\,x_i^{-k_i}\\\times
(a^{r-1}q^{1-|\bk|}/E)_{|\bk|}
\left(\frac d{a^{r-1}}\right)^{|\bk|}\,
q^{(r+1)\binom{|\bk|}2-\sum_{i=1}^r\binom{|\bk|-k_i}2}\\
=\prod_{i,j=1}^r\frac {(ax_iq^{1+k_i}/e_jx_j)_\infty}{(ax_iq/e_jx_j)_\infty}
\prod_{i=1}^r\frac{(a^{r-1}e_ix_iq^{1-|\bk|}/dE)_\infty\,(a^rx_iq/dE)_\infty}
{(a^{r-1}e_ix_iq/dE)_\infty\,(a^rx_iq^{1+k_i-|\bk|}/dE)_\infty}\,
x_i^{-k_i}\\\times
\frac{(a^{r-1}q^{1-|\bk|}/E)_\infty}{(a^{r-1}q/E)_\infty}
\left(\frac d{a^{r-1}}\right)^{|\bk|}\,
q^{\binom{|\bk|}2+|\bk|^2-\sum_{i=1}^r\binom{k_i+1}2},
\end{multline}
followed by multiple applications of the $q$-binomial theorem.
The last expression thus appears implicitly as a factor of the
summand of the  multilateral series over $k_1,\dots,k_r$. 
However, it is most important that the ``quadratic'' powers of $q$
in the last line of \eqref{auxid}, specifically
\begin{equation*}
q^{\binom{|\bk|}2+|\bk|^2-\sum_{i=1}^r\binom{k_i+1}2}=
q^{e_2(\bk)+2\binom{|\bk|}2}
\end{equation*}
(where $e_2(\bk)$ is the second elementary symmetric function of $\bk$),
ensure absolut convergence of the multiple power series
in $1/e_1,\dots,1/e_r$.
\end{proof}

\end{document}